\documentstyle[12pt]{article}

\newcommand{\const}{\mathop{\rm const}\limits}

\newcommand{\supp}{\mathop{\rm supp}\limits}

\begin{document}

\begin{center}

{\bf FUNDAMENTAL FUNCTION FOR }\\

\vspace{3mm}

{\bf GRAND LEBESGUE SPACES. } \\

\vspace{4mm}

                E.Ostrovsky, L.Sirota  \\

\vspace{3mm}

 \ Department of  Mathematic, Bar-Ilan University,  Ramat Gan, 52900, Israel, \\
 e-mails:  \ eugostrovsky@list.ru, \ sirota3@bezeqint.net \\

\vspace{3mm}
                    {\sc Abstract.}\\

 \end{center}

 \vspace{3mm}

 \ We investigate in this short article the fundamental function for the so-called Grand Lebesgue Spaces (GLS)
and show in particular a one-to-one and mutually continuous accordance between its fundamental and generating function.

  \vspace{3mm}


\vspace{3mm}


\vspace{2mm}

\section{Notations. Statement of problem.}

\vspace{4mm}

 \ {\bf A.} A triplet $ (X, \cal{B}, \mu ),  $  where $ X = \{x\}  $ is arbitrary set,
$ \cal{B} $ is non-trivial certain sigma-algebra of subsets $ X  $ and  $ \mu  $ is probabilistic: $ \mu(X) = 1 $
 diffuse non-negative completely additive measure defined on the $  \cal{B}. $ \par

  \ The non-probabilistic case $  \mu(X) = \infty  $ will be consider further. \par

 \ Recall that the measure $  \mu $ is said to be diffuse iff for arbitrary measurable set $  A_1 \in \cal{B} $
with positive measure:  $ \ \mu(A_1) > 0  $
there exists it subset $  A_2 \subset A_1 $ such that $ \mu(A_2) = \mu(A_1)/2.  $ \par

 \ We denote as usually for any arbitrary measurable function $ f: X \to R   $

 $$
 |f|_p =  \left[ \int_{ X } f(x)|^p \ \mu(d x)  \right]^{1/p}, \ p \ge 1;
 $$

$$
L_p = \{ f, \ |f|_p < \infty  \}.
$$

\vspace{3mm}

 \ {\bf B.} The so-called Grand Lebesgue Space (GLS) $ G \psi $  with norm $  ||\cdot ||G\psi   $ is defined
 ( not only in this article) as follows:

$$
G \psi = \{f, \ ||f||G \psi < \infty \}, \  ||f||G \psi \stackrel{def}{=} \sup_{ p \ge 1 } \left[ \frac{|f|_p}{\psi(p)} \right].
\eqno(1.1)
$$
 \ Here $ \psi = \psi(p), \ 1 \le p < \infty $ is some continuous strictly increasing function such that
 $ \lim_{p \to \infty} \psi(p) = \infty. $\par \par

 \ The detail investigation of this spaces (and more general spaces) see in \cite{Liflyand1},  \cite{Ostrovsky8}. See also
 \cite{Fiorenza1}, \cite{Fiorenza2}, \cite{Iwaniec1}, \cite{Iwaniec2}, \cite{Kozachenko1} etc. \par

\vspace{4mm}

 \ The case when in (1.1) supremum is calculated over {\it finite } interval is investigated in
 \cite{Liflyand1}, \cite{Ostrovsky8}, \cite{Ostrovsky9}:

$$
G_b \psi = \{f, \ ||f||G_b \psi < \infty \}, \  ||f||G_b \psi \stackrel{def}{=}
\sup_{ 1 \le p < b } \left[ \frac{|f|_p}{\psi(p)} \right],
 \  b = \const > 1,  \eqno(1.2)
$$
but in (1.2) $ \psi = \psi(p) $ is  continuous  function in the semi-open interval $ 1 \le p < b $
such that $ \inf_{p \in (1,b)} \psi(p) > 0. $ \par

 \ We will denote

$$
[1,b) := \supp \psi(\cdot),
$$
 or simple $ b  = b(\psi):= \supp \psi(\cdot),  $ including the case $  b = \infty. $ \par

\vspace{4mm}

 \ {\bf Definition 1.1.} The function $  = \psi(p)  $ which appeared in (1.1) and (1.2), will be named as
 {\it  generating function } for the correspondent Banach space $  G\psi.  $\par

\vspace{4mm}

 \ An used further example:

 $$
 \psi^{(\beta,b)}(p) = (b-p)^{-\beta}, \ 1 \le p < b,  \beta = \const > 0; \ b = \const > 1,
 $$

 $$
  G_{\beta,b}(p) := G_b\psi^{(\beta,b)}(p).
 $$

\vspace{3mm}

 \ {\bf C.} We denote as ordinary  for any measurable set $ A, \ A \in \cal{B} $ it indicator function by
$ I(A) = I_A(\omega).   $\par

\vspace{3mm}

 \ {\bf D.}  The Grand Lebesgue Spaces   $ \{ G \psi \} $  are rearrangement invariant in the classical definition, see
e.g. \cite{Bennet1}, chapter 1. Therefore, its fundamental function $ \phi_{G(\psi) }(\delta), \ \delta \ge 0 $ is correctly
defined  in the considered case as follows:

$$
\phi_{G(\psi) }(\delta) \stackrel{def}{=} \sup_{p \in \supp \psi} \left[ \frac{\delta^{1/p}}{\psi(p)} \right], \eqno(1.3)
$$
see \cite{Bennet1}, chapters 2 and 5. \par

 \ For instance,

$$
\phi_{G(\psi)}(1) =   \frac{1}{\inf_{p \in \supp \psi} \psi(p)}. \eqno(1.3a)
$$

 \ Note also

$$
\phi_{G(C \cdot \psi)} = \phi_{G(\psi)}/C, \ C = \const >0. \eqno(1.3b)
$$

\vspace{3mm}

 \ This notion play a very important role in the functional analysis, \cite{Bennet1}, \cite{Ostrovsky11}, \cite{Ostrovsky12};
in the theory of interpolation of operators, \cite{Bennet1}, \cite{Fiorenza1}, \cite{Fiorenza3}, in the theory of probability
\cite{Talagrand1},  \cite{Ostrovsky2}, \cite{Ostrovsky3}, \cite{Ostrovsky5}; in the theory of Partial Differential equations
\cite{Fiorenza3}, \cite{Iwaniec2}; in the theory of martingales \cite{Ostrovsky9};
in the theory of approximation, in the theory of random processes etc. \par

\vspace{3mm}

 \ {\bf E.} Let $ g = g(p), \ p \in (a,b), 1 \le a < b \le \infty  $ be some numerical valued continuous strictly increasing
(or decreasing)  function. The {\it inverse }  function will be denoted by $  g^{(-1)}(z), \ g(a) \le z \le g(b),  $ in
contradistinction to the usually  notation $  g^{-1}(p) = 1/g(p). $ \par

\vspace{3mm}

 \ {\bf F.} The Young-Fenchel, or Legendre transform $   g^*(q)  $ for the  function $ g = g(p)  $ one can to define

$$
 g^*(q) \stackrel{def}{=} \sup_{p \in \supp g} (p \ |q| - g(p)). \eqno(1.4)
$$

\vspace{4mm}

 \ {\bf Our goal in this short report is to establish a one-to-one and mutually continuous connection between  the
 fundamental and generating functions for the Grand Lebesgue Spaces.}\par

\vspace{4mm}

 \ In some previous works: \cite{Fiorenza2}, \cite{Kufner1}, chapter 8; \cite{Liflyand1},
\cite{Ostrovsky8}, \cite{Ostrovsky6}, \cite{Ostrovsky11}  these function  was evaluated and applied in many practical cases.\par

\vspace{3mm}

\section{Main result.}

\vspace{3mm}

 \ {\bf Problem A.  Let the generating function $  \psi $ be a given: $  \psi \in G\Psi(a,b), 1 \le a < b \le \infty.    $
 Find the fundamental function for the correspondent Grand Lebesgue Space $ G\psi. $ } \par

\vspace{3mm}

 \ Suppose the function

$$
p \to \frac{p}{\psi(p)}, \ p \in (a,b)
$$
is strictly increasing; and define therefore the function

$$
\nu(p) = \nu_{\psi}(p) = \left[ \frac{p}{\psi(p)}   \right]^{(-1)}, \ p \in \supp \psi, \eqno(2.1)
$$
and  $ \nu(p) = + \infty  $  otherwise.\par

 \ Introduce also the following Young-Orlicz function

$$
N(u) = N_{\psi}(u) := \exp \left(  \nu_{\psi}^*(u)  \right) - \exp \left(  \nu_{\psi}^*(0)  \right), \eqno(2.2)
$$
and define finally

$$
\theta(\delta) = \theta_{\psi}(\delta) \stackrel{def}{=}  \frac{1}{N_{\psi}^{(-1)}(1/\delta)}, \ \delta > 0. \eqno(2.3)
$$

\vspace{3mm}

{\bf Proposition 2.1.} {\it We propose under formulated above conditions, for instance, $  \mu(X) = 1, $
diffuseness of the measure $  \mu, $  and in the case when  } $  b = \infty $

\vspace{3mm}

$$
\phi_{G(\psi)}(\delta) = \theta_{\psi}(\delta), \ \delta > 0. \eqno(2.4)
$$

\vspace{3mm}

\ {\bf Remark 2.1.} The equality (2.4) is more convenient than source definition (1.3). In particular,
it allows for a relatively simple inversion. \par

\vspace{3mm}

 \ {\bf Proof} is very simple; it based on the computation of the fundamental function for Orlicz spaces, see
the book of  Krasnosel'skii M.A. and Rutickii Ya.B.  \cite{Krasnosel'skii1}, chapter 3; see also the classical
monographs \cite{Rao1},  \cite{Rao2}. \par

 \ In detail, it is proved in particular in the articles \cite{Liflyand1}, \cite{Ostrovsky6}, \cite{Ostrovsky8}
that under our conditions the Grand Lebesgue Space $  G\psi $ coincides with certain  Orlicz space over source probability
triplet $ (X, \cal{B}, \mu ) $  relative the Young-Orlicz function $  N_{\psi}(u). $ \par

 \ We  deduce reducing considered case to the well-known calculation of fundamental function for Orlicz space,
\cite{Krasnosel'skii1}, chapter 3

$$
\phi_{G(\psi)}(\delta) = \frac{1}{N^{(-1)}_{\psi}(1/\delta)} = \theta_{\psi}(\delta), \ \delta > 0, \eqno(2.4)
$$
Q.E.D. \par

\vspace{4mm}

\ {\bf An inverse problem B. \ Let the fundamental function $  \phi_{G\psi}(\delta) = \phi(\delta) $
 be a given. Find the correspondent generating function $  \psi(p). $ } \par

 \ A first restrictions: the function $  \phi = \phi(\delta) $  is strictly increasing and continuous; in particular
$  \phi(0+) = \phi(0) = 0.  $ \par

 \ We find from the equality (2.4)

$$
N_{\psi}(1/\delta) = \left( \frac{1}{\phi(\delta)} \right)^{(-1)}, \eqno(2.5)
$$
or equivalently

$$
N_{\psi}(z) = \left( \frac{1}{\phi(\delta)} \right)^{(-1)} \left/_{ \delta := 1/z} \right/. \eqno(2.5a)
$$

 \vspace{3mm}

 \ A second restriction: the function

 $$
 V^*(z) = \ln (C + N_{\psi}(z) ), \ z \ge 0, \eqno(2.6)
 $$
where $ V^*(0) = \ln C,  $ is continuous and upward convex.\par
 \ It follows immediately from (2.6)  by virtue of theorem of Fenchel-Moraux

$$
 V(z) =  \left\{\ln (C + N_{\psi}(z) ) \right\}^*, \ z \ge 0. \eqno(2.7)
 $$

 \ Since

$$
V(p) = \left[  \frac{p}{\psi(p)}  \right]^{(-1)},
$$
we  derive finally \\

{\bf Proposition 2.2.} We conclude under formulated in this pilcrow conditions

$$
\psi(p) = \frac{p}{V^{(-1)}(p)}. \eqno(2.8)
$$

\vspace{3mm}

\section{The case of infinite measure.}

\vspace{3mm}

 \ The case  when $  \mu(X) = \infty $  is more complicated.\par

 \ Recall first of all definition and some facts about the so-called Exponential Orlicz Spaces (EOS), see
for example \cite{Ostrovsky10}. \par

 \ Let $ N = N(u) $ be an
$ N \ - $  Young-Orlicz's function, i.e. downward convex, even, continuous function
differentiable for all sufficiently greatest values $ u, \ u \ge u_0, \ u_0 = \const > 0, $
strongly  increasing along the right semi-axis
and such that $ N(u) = 0 \ \Leftrightarrow  u = 0; \ u \to \infty \ \ \Rightarrow
dN(u)/du \ \to \ \infty. $ We can say that  $ N(\cdot) $ is an exponential
Orlicz function,
briefly, $ N(\cdot) \in \ EOF, $ if $ N(u) $ has a form of a continuous
differentiable strongly increasing downward convex function $ W = W(u) $ in the domain $ [2,\infty] $
such that $ u \to \infty \ \Rightarrow W^/(u) \to \infty $ and

$$
N(u) = N(W,u) = \exp(W(\log |u|)), \ |u| \ge e^2.
$$
 For the values $ u \in [-e^2,e^2] $ we define $ N(W,u) $ arbitrarily,  but so that the function
$ N(W,u) $ is even, continuous, convex, strictly increasing along the right semi-axis
and so that $ N(u) = 0 \
\Leftrightarrow u = 0. $ We denote the correspondent Orlicz space on $ (X,\mu) $ with a
 measure $  \mu  $ and with $ N \ - $ function of the form $ N(W,u) $ as
 $ L(N) = EOS(W); \ EOS = \cup_{W} \{EOS(W)\} $ (exponential Orlicz's space).\par
 \ For example, let $ m = \const > 0, \ r = \const \in R^1, $
$$
N_{m,r}(u) = \exp \left[|u|^m \ \left(\log^{-mr}(C_1(r) + |u| ) \right) \right] -1,
$$
$ C_1(r) = e, \ r \le 0; \ C_1(r) = \exp(r), \ r > 0. $ Then $ N_{m,r}(\cdot) \in EOS. $
In the case $ r = 0 $ we can write $ N_m = N_{m,0}. $ \par
 Recall that the Orlicz's norm on the arbitrarily measurable space $ (X,A,\mu) $
$ ||f||L(N) = ||f||L(N,X,\mu) $ can be calculated by the following formula (see, for example, \cite{Krasnosel'skii1}, p. 66;
\cite{Rao1}, p. 73 )
$$
||f||L(N) = \inf_{v > 0}  \left \{ v^{-1} \left(1 + \int_X N( v|f(x)|) \ \mu(dx) \right) \right \}.
$$

 \ Let $ \alpha $ be arbitrary number, $ \alpha = \const \ge 1, $ and $ N(\cdot)
\in EOS(W) $ for some $ W = W(\cdot). $  For such a function
$ N = N(W,u) $ we denote by $ N^{(\alpha)}(W;u) = N^{(\alpha)}(u) $ a new
 Young-Orlicz's  function $ N^{(\alpha)}(u) $ such that

$$
N^{(\alpha)}(u) = C_1 \ |u|^{\alpha}, \ \ |u| \in [0, C_2];
$$
$$
N^{(\alpha)}(u) = C_3 + C_4 |u|, \ \ |u| \in (C_2, C_5];
$$
$$
N^{(\alpha) }(u) = N(u), \ \ |u| > C_5, \ \ 0 < C_2 < C_5 < \infty,\eqno(3.1)
$$
$$
C_{1,2,3,4,5} = C_{1,2,3,4,5}(\alpha,N(\cdot)).
$$
 \ In the case of $ \alpha = m(j+1), \ m > 0, \ j = 0,1,2, \ldots $ the function $ N^{(\alpha)}_m(u) $
is equivalent to the following Trudinger's function:
$$
N_m^{(\alpha)}(u) \sim N_{[m]}^{(\alpha)}(u) = \exp \left(|u|^m \right) -\sum_{l=0}^{j}
u^{ml}/l!.
$$
 \ This method is described in {\cite{Taylor1},
 p. 42-47. These Orlicz spaces are applicable to the theory of non-linear partial differential equations.\par

\vspace{3mm}

  \ We denote  hereinafter generally by $  C_k = C_k(\cdot), k = 1,2,\ldots $ some positive finite
essentially constructive constants, and by $ C, C_0 $ non-essentially  constants, also constructive.
 \ We  proved the existence of constants $ C_{1,2,3,4,5} = C_{1,2,3,4,5}(\alpha,
N(\cdot)) $ such that $ N^{(\alpha)} $ is a new exponential $ N \ $
Orlicz's function in \cite{Ostrovsky10}. We denote classical absolute constants by the symbols $ K_j $.\par

 \  Now we introduce some {\it new } Grand Lebesgue Spaces. Let $ \psi = \psi(p),
\ p \ge \alpha, \alpha = \const \ge 1 $
be a continuous positive $ \psi(\alpha) > 0 $  finite strictly increasing function such
that the function $ p \to p \log \psi(p) $ is downward convex, and
$$
\lim_{p \to \infty} \psi(p) = \infty.
$$
We denote the set of all these functions  by $ \Psi; \ \Psi = \{\psi\}. $ A particular case

$$
\psi(p) = \psi(W;p) = \exp(W^*(p)/p),
$$
where
$$
W^*(p) = \sup_{z \ge \alpha} (pz - W(z))
$$
is so-called Young-Fenchel, or Legendre transform of $ W(\cdot). $ It follows from the
theorem of Fenchel-Moraux that in this case
$$
W(p) = \left[ p \ \log \psi(W;p) \right]^*, \ \ p \ge p_0 = const \ge 2,
$$
and, consequently,  for all $ \psi(\cdot) \in \Psi $ we introduce a correspondent Young-Orlicz $ N \ - $
function by the equality:
$$
N([\psi]) = N([\psi],u) =
\exp \left \{\left[p \log \psi(p) \right]^*(\log u) \right \}, \ u \ge e^2.
$$

\vspace{3mm}

{\bf Definition 3.1. } We introduce for such arbitrary function $ \psi(\cdot) \in \Psi $
the so-called $ G(\alpha;\psi) \ $  norms  and correspondent Banach GLS space $ G(\alpha; \psi) $ as a
set of all measurable (complex) functions with finite norms:

$$
||f||G(\alpha; \psi) = \sup_{p \ge \alpha} (|f|_p/\psi(p)). \eqno(3.2)
$$

 \ For instance, $ \psi(p) $ may be $ \psi(p) = \psi_m(p) = p^{1/m}, \ m = \const > 0; $ in
this case, we can write $ G(\alpha, \psi_m) = G(\alpha,m) $ and
$$
||f||G(\alpha,m) = \sup_{p \ge \alpha}  \left(|f|_p \ p^{-1/m} \right).
$$

\vspace{4mm}

{\bf Theorem A, see } \cite{Ostrovsky10}.  {\it  Let the measure $  \mu $ be diffuse,
 $  \mu(X)  = \infty, \alpha = \const \ge 1  $ and $ \psi \in \Psi. $
We assert that the norms Orlicz-Luxemburg norm $ ||\cdot|| L(N^{(\alpha)}, [\psi]) $ and Grand Lebesgue Space norm
$||\cdot|| G(\alpha, \psi), \ \alpha \ge 1 $ are equivalent.} \par

\vspace{3mm}

 \ Arguing similarly to the second section, we obtain the following  result.

\vspace{3mm}

{\bf Proposition 3.1.} {\it We propose under  conditions of theorem A }

\vspace{3mm}

$$
\phi_{G(\psi)}(\delta) \asymp \theta_{\psi}(\delta), \ \delta > 0.
$$

\vspace{3mm}

 {\bf Remark 3.1.} Note that this case $ \delta \in (0, \infty),  $ in contradiction to the proposition 2.1, where
it in naturally to take  $ \delta \in (0,1). $ \par

\vspace{3mm}

\section{Concluding remarks. Open problems.}

\vspace{3mm}

 \ It is interest by our opinion to investigate the notion of fundamental function and also its
relation with generating function  for the so-called mixed, or equally anisotropic Grand Lebesgue Spaces.  \par
 \ Recall that the definition of mixed, or equivalently anisotropic ordinary Lebesgue Spaces appeared at first
in the article \cite{Benedek1}  and was investigated in detail in the classical books   \cite{Besov1},  \cite{Besov2}. \par

  \ The anisotropic Grand Lebesgue Spaces as a slight generalization of $ L_p $ spaces arises in turn in \cite{Ostrovsky12}
with the correspondent fundamental function;  in the preprint   \cite{Ostrovsky14}  both this notions
 was applied in the operator's theory. \par

\end{document}